\documentclass[12pt,a4paper,twoside]{amsart}
\usepackage{amsfonts, amsthm, amsmath, amssymb}
\usepackage{mathrsfs,amsmath}
\usepackage{hyperref}
\usepackage{bbm}
\hypersetup{colorlinks=false}

\usepackage[margin=2.9cm]{geometry}

\usepackage{helvet}

\usepackage{amsthm}
\usepackage{lineno}
\newtheorem{theorem}{Theorem}
\newtheorem{lemma}{Lemma}[subsection]
\newtheorem{remark}{Remark}
\newtheorem{application}{Application}

\newtheorem{corollary}{Corollary}

\DeclareMathOperator{\lcm}{lcm}

\begin{document}

\author{ Aritra Ghosh }
\title{Subconvex bound for Rankin-Selberg $L$-functions in prime power level}

\address{Aritra Ghosh \newline  Alfr\'ed R\'enyi Institute of Mathematics, Budapest, Re\'altanoda utca 13-15, 1053, Hungary; email: aritrajp30@gmail.com}

\begin{abstract}
Let $f$ be a $p$-primitive cusp form of level $p^{4r}$, where local representation of $f$ be supercuspidal at $p$, $p$ being an odd prime, $r\geq 1$ and $g$ be a Hecke-Maass or holomorphic primitive cusp form for $\mathrm{SL}(2,\mathbb{Z})$. A subconvex bound for the central values of the Rankin-Selberg $L$-functions $L(s, f \otimes g )$ is given by $$ L (\frac{1}{2}, f \otimes g ) \ll_{g,\epsilon}p^{\frac{23r}{12} +\epsilon}.$$

\end{abstract}

\maketitle

\section{Introduction}
In this paper, we are interested in the program – initiated by various mathematicians to solve the subconvexity problem for Rankin–Selberg $L$-functions (see \cite{HM}, \cite{KMV}, \cite{LLY}, \cite{PM}, \cite{MV}, \cite{PS}) in the level aspect. We continue that program here, in this paper also, using the delta method which was introduced and developed by Munshi in his papers (see \cite{Mun3}, \cite{Mun}, \cite{Mun4}).

Let $f$ be a $p$-primitive cusp form which is supercuspidal at $p$ with level $p^{4r}$, $p$ being an odd prime, $r\geq 1$ and $g$ be a Hecke-Maass or holomorphic primitive cusp form for $\mathrm{SL}(2,\mathbb{Z})$. Then the $L$-function associated with $f \otimes g $ is given by
$$L(s, f \otimes g ) =  \zeta_{p}(2s)\sum_{n=1}^{\infty}\frac{\lambda_{f}(n)\lambda_{g}(n)}{n^s}, \, \, \, \, \zeta_{p}(s)= \prod_{q\neq p }\left( 1-\frac{1}{q^s}\right)^{-1},$$
for $\mathrm{Re}(s) > 1$, which can be analytically extended to the whole complex plane $\mathbb{C}$ and satisfies a functional equation relating $s$ and $1-s$ and for the last product in the previous line, $q$ runs over all primes except $p$. In this context the convexity bound is $$L(\frac{1}{2}, f \otimes g ) \ll_{g,\epsilon}p^{2r +\epsilon},$$ for any $\epsilon >0$ which can be obtained by using the approximate functional equation and the Phragm\'en-Lindel\"of principle. Then the subconvexity bound we are seeking is a bound of the form
$$L(\frac{1}{2}, f \otimes g ) \ll_{g,\epsilon}p^{2r-\delta +\epsilon},$$

\noindent
with some absolute $\delta>0$, the implied constant depending
on $g$ and the parameter at infinity of $f$ (i.e., the weight or the Laplacian eigenvalue but not on the level of $f$). This problem was first addressed by Kowalski-Michel-Vanderkam (see \cite{KMV}) under the following assumptions:
\begin{itemize}
    \item $f$ is a holomorphic cusp form;
    \item the conductor $q^{*}$ of $\chi_{f}\chi_{g}$ is at most $q^{\frac{1}{2}-\eta}$ for some $\eta>0$, the corresponding exponent $\delta$ then depending on $\eta$ where $\chi_{f}$ is the nebentypus of $f$ and $\chi_{g}$ is the nebentypus of $g$.
\end{itemize}

Then the second condition was essentially removed in \cite{PM} under the assumptions:
\begin{itemize}
    \item $g$ is a holomorphic cusp form;
    \item $\chi_{f}\chi_{g}$ is non-trivial.
\end{itemize}

Later Harcos-Michel (see \cite{HM}) considered the previous problem after removing the first assumption that $g$ is a holomorphic cusp form though $\chi_{f}\chi_{g}$ being non-trivial. Michel-Venkatesh have also considered this problem in their breakthrough work \cite{MV}. If the level of the form $f$ is a prime then this problem has been considered in \cite{CR}, \cite{GR}. Raju (see \cite{CR}) did the problem by taking $g$ to be cuspidal form though his method does not work when $g$ is an Eisenstein series. For some time it was not clear whether the $\delta$-method approach will work in the case of $g$ non-cuspidal. Recently Aggarwal, Kumar, Kwan, Leung, Li and Young (see \cite{AKKLLY}) solved this subconvexity problem in level aspect using the $\delta$-method and here they can take $g$ to be an Eisenstein series. The author has also considered this problem as a $\mathrm{GL}(2)\times \mathrm{GL}(2) \times \mathrm{GL}(1)$ problem in \cite{AG1}. Also in this context one can see the work of Holowinsky-Munshi (see \cite{RHM}), Holowinsky-Templier (see \cite{HT}), Hou-Zhang (see \cite{HZ}) for the hybrid level aspect case. 

In this paper we will consider the case when the level of the form $f$ is a prime power, getting the following result:

\begin{theorem}\label{mt}
     Let $f$ be a $p$-primitive cusp form of level $p^{4r}$ with the local representation of $f$ be supercuspidal at $p$, $p$ being a sufficiently large odd prime, $r\geq 1$ and $g$ be a Hecke-Maass or holomorphic primitive cusp form for $\mathrm{SL}(2,\mathbb{Z})$. Then we have
    $$L(\frac{1}{2}, f \otimes g ) \ll_{g,\epsilon}p^{\frac{23r}{12}+\epsilon}.$$
\end{theorem}

\begin{remark}
 Here we are considering $f,g$ to be holomorphic cusp forms for simplicity as we can proceed similarly for Hecke-Maass cusp forms. Also here we will consider the case for $r=1$ as for $r>1$, we can proceed similarly. For $r\geq 2$ one can use \cite[Theorem $1$]{FS} and \cite[Remark $2$]{FS} to get the square root cancellation for the character sum \eqref{a} in the off-diagonal case. If the form $f$ has level $p^u$ for some $u \in \mathbb{N}$, then for the case $u\geq 4$ and $4\nmid u$ the author hopes that one can do address the problem for general case following the same method and doing conductor lowering by taking the conductor of size $p^{[\frac{u}{4}]}$. The author also plans to address the general problem in a separate article in a different approach.
\end{remark}

 

\begin{remark}
    Our result gives an improvement of the results done in \cite[Theorem $1.1$]{KMV}. The main input of our work is to use the Jutila's circle method (see \ref{sec3}) in a certain way so that after an application of the Poisson summation formula, only zero frequencies survive. By repeated application of the Cauchy-Schwarz's inequality and change of variables, we separate the oscillations and get new character sums.
\end{remark}

\subsection*{Notation}
In this article, by $A \ll B$ we mean $|A| \leq C |B|$ for some absolute constant $C > 0$, depending only on $g,\epsilon$ and notation `$X\asymp Y$' will mean that $Y p^{-\epsilon}\leq X \leq Y p^{\epsilon}$. 

\begin{application}
As an application of our theorem \eqref{mt}, we improve the bound obtained in \cite{KMV} for the problem of distinguishing modular forms based
on their first Fourier coefficients.
\begin{corollary}
From the proof of the corollary $1.3$ of \cite{KMV}, we have, for $f$ as given in our theorem \eqref{mt}. Then there exists a constant $C = C( f, \epsilon)$ such that for any primitive cuspidal newform $g$ for $\mathrm{SL}(2,\mathbb{Z})$ and for any primitive Dirichlet character $\chi(n)$ having modulus $N$, odd, there exists $n \leq C \, N^{1-\frac{1}{24} +\epsilon}$ with $(n,N)=1$, such that
$$\lambda_{f}(n)\neq \lambda_{g}(n) .$$
\end{corollary}
\begin{proof}
    The proof is identical to the proof of Corollary $1.3$ \cite{KMV}. Use Theorem \eqref{mt} in place of $[$\cite{KMV}, Theorem $1.1$$]$ in the proof.
\end{proof}
\end{application}

\subsection*{Acknowledgement}
The author thanks Ritabrata Munshi and Gergely Harcos for their helpful comments and suggestions. The author is also grateful to the Alfr'ed R'enyi Institute of Mathematics, Budapest, for the excellent research environment. The author thanks Debargha Banerjee for his help with the root number computations in the Subsection \ref{root}. The author is also grateful to Mayukh Dasaratharaman for his helpful comments and suggestions. 

\section{Preliminaries}\label{sec2}  
\subsection{Preliminaries on holomorphic cusp forms}
  Let $f$ be a $p$-primitive cusp form where the level of $f$ be $p^{4}$ with the local representation of $f$ being supercuspidal in $p$, $p$ being an odd prime, and $g$ be a Hecke-Maass or holomorphic primitive cusp form for $\mathrm{SL}(2,\mathbb{Z})$ 
 \subsubsection{Voronoi summation formula.}  We will use the following Voronoi summation formula. In the case of a full level (see \cite[Proposition $2$]{BH}) we have

 \begin{lemma}\label{V} Let $a$ and $c$ be coprime positive integers, and let $F:(0,\infty)\mapsto \mathbb{C}$ be a
smooth function of compact support. Then
 \begin{equation}\label{eq9}
 c\sum_{n=1}^{\infty}\lambda_{g}(n)e\left( \frac{an}{c}\right)F(n)=\sum_{n=1}^{\infty}\lambda_{g}(n)\sum_{\pm}e\left(\mp \frac{\overline{a}n}{c}\right)\int_{0}^{\infty}F(x)J_{g}^{\pm}\left( \frac{4\pi \sqrt{\sqrt{nx}}}{c}\right)\mathrm{d}x,
 \end{equation}

 \noindent
 where 

$$J_{g}^{+}(x):= 2\pi i^k J_{k-1}(x), \, \, \, \, J_{g}^{-}(x):= 0,$$

\noindent
if $g$ is a holomorphic cusp form of level $1$ and weight $k$;

$$J_{g}^{+}(x):= \frac{-\pi}{\cosh (\pi r)}(Y_{2ir}(x)+ Y_{-2ir}(x)), \, \, \, \, J_{g}^{-}(x):= \varepsilon_{g}4\cosh (\pi r)K_{2ir}(x),$$

\noindent
if $g$ is a Maass cusp form of level $1$, weight $0$, Laplacian eigenvalue $\frac{1}{4}+r^2$ and sign $\varepsilon_{g}\in \{ \pm 1\}$. For $g=E_{r} \, \, (r \in \mathbb{R})$ the same formula holds with $J_{g}^{\pm}$ as in the Maass case (with $\varepsilon_{g}=1$), except that on the right-hand side the following polar term has to be added:

$$\sum_{\pm}\zeta (1\pm 2ir)\int_{0}^{\infty}\left( \frac{x}{c^2}\right)^{\pm ir}F(x) \mathrm{d}x \, \, \, \, \textit{ for } r\neq 0,$$

$$\int_{0}^{\infty}\left( \log \left( \frac{x}{c^2}\right)+2\gamma\right) F(x)\mathrm{d}x \, \, \, \, \textit{ for } r=0.$$
 
\end{lemma}

\subsubsection{Bounds for the Bessel functions}\label{bessel} For $x\gg 1$, we have
$$J_{k-1}(2\pi x)= W_{k}(x)\frac{e(x)}{\sqrt{x}}+\overline{W}_{k}(x)\frac{e(-x)}{\sqrt{x}},$$

\noindent
for some function $W_{k}$ satisfying $x^j (\partial^j W_{k})(x)\ll_{j,k}1$.

\subsubsection{The Ramanujan bound on average}\label{RAS} Here we note the following Ramanujan bound on average (see \cite{IW}):
\begin{equation}\label{RA}
    \sum_{n\ll N}|\lambda_{f}(n)|^2 \ll p^\epsilon N^{1+\epsilon}.
\end{equation}

 \subsection{Circle Method}\label{circle}
 
 Here in this paper we shall use Jutila's circle method (see \cite{PA}, \cite{MZ}). For any set $S \subset R$, let $I_{S}$ denote the associated characteristic function, i.e. $I_{S}(x)=1$ for $x \in S$ and $0$ otherwise. For any collection of positive integers $\Phi \subset [Q,2Q]$ (which we call the set of moduli), where $Q \geq 1 $ and a positive real number $\delta$ in the range $Q^{-2} \ll \delta \ll Q^{-1}$, we define the function
$$\tilde{I}_{\Phi , \delta}(x):= \frac{1}{2\delta L}\sum_{q \in \Phi} \, \, \sideset{}{^*}\sum_{d \bmod q}I_{[\frac{d}{q}-\delta , \frac{d}{q}+ \delta]}(x),$$
where $I_{[\frac{d}{q}-\delta , \frac{d}{q}+ \delta]}$ is the indicator function of the interval $[\frac{d}{q}-\delta , \frac{d}{q}+ \delta]$. Here $L := \mathop\sum_{q \in \Phi}\phi (q)$ (then roughly we have $L \asymp Q^2$) and we will choose $\Phi$ in such a way that $L \asymp Q^{2-\epsilon}$.
 
Then this becomes an approximation of $I_{[0,1]}$ in the following sense:
\begin{lemma}\label{2.2} We have
$$\int_{\mathbb{R}}\left| I_{[0,1]}(x) - \tilde{I}_{\Phi , \delta}(x)\right|^2 dx \ll \frac{Q^{2+\epsilon}}{\delta L^2},$$
where $I$ is the indicator function of $[0,1]$.
\end{lemma}

\noindent
This is a consequence of the Parseval theorem from Fourier analysis (see \cite{PA}, \cite{Mun3}).

\subsection{Root number}\label{root} Let $f$ be a newform of level $p^j$ ($j>2$) with trivial nebentypus. We denote by $\pi_{f}$ the automorphic representation (attached to $f$) of the ad\'ele group $\mathrm{GL}_{2}(\mathbb{A}_{\mathbb{Q}})$. It is well-known that $\pi_{f}$ decomposes as a restricted tensor product $\pi_{f}=\otimes_{v}^{\prime}\pi_{f,v}$ where $v$ runs over
all places of $\mathbb{Q}$. Here $\pi_{f,v}$ is an irreducible admissible representation of $\mathrm{GL}_{2}(\mathbb{Q}_{v})$ and $\mathbb{Q}_{v}$ denotes the completion of $\mathbb{Q}$ at $v$.

Let $\chi$ be a Dirichlet character modulo $p^r$ with $r<j$. We assume that $r$ is quite small compared
to $j$ ($r\leq \left[ \frac{j}{4} \right] $) and $p$ is a large prime number. The question we are interested in is the following: What is the epsilon factor/root number of $f\otimes \chi$? By class field theory, the character $\chi$ can be identified with a character of the id\'ele group, that is, characters $\{\chi_{q} \}_{q}$ with $\chi_{q}:\mathbb{Q}_{q}^{\times}\rightarrow \mathbb{C}^{\times}$, see \cite[p. 195]{KL} for more details. Now to compute the local epsilon factors of $f$ at each prime $p$, we choose an additive character $\phi$
of $\mathbb{Q}_{p}$ of conductor $0$. Note that the global epsilon factor of $f$ is a product of all local epsilon factors and that does not depend on the additive character chosen. Then one can show that, for the prime $p$,
    \begin{equation*}
    \varepsilon (\pi_{f,p}\otimes \chi_{p},\phi)= 
\begin{cases}
    \chi_{p}(p^j u)\varepsilon (\pi_{f,p},\phi),& \text{ if $p$ is an unramified supercuspidal prime}, \\
    \chi_{p}(p^{2j}u)\varepsilon (\pi_{f,p},\phi),& \text{ if $p$ is an ramified supercuspidal prime} ,
\end{cases}
    \end{equation*}

\noindent
for some $u \in \mathbb{Z}_{p}^{\times}$. Also for other primes $q\neq p$, we have

    $$\varepsilon (\pi_{f,q}\otimes \chi_{q},\phi)=\varepsilon (\pi_{f,q},\phi).$$

\noindent
Hence we have

\begin{equation}\label{ep}
    \varepsilon (\pi_{f}\otimes \chi)= \chi(\mathfrak{a})\varepsilon (\pi_{f}) ,
\end{equation}

\noindent
for some constant $p \nmid \mathfrak{a}$, depending on $f$ only. Without loss of generality, or our computation we will take $\mathfrak{a} =1$ which can be seen from the computation.

 \section{Setting-up the circle method :}\label{sec3}
 \subsection{Approximate functional equation}By means of approximate functional equation and dyadic subdivision, we have 
\begin{equation}
   L(\frac{1}{2}, f \otimes g ) \ll_{\epsilon, A} \sup_{N \ll p^{2 + \epsilon}}\frac{S(N)}{\sqrt{N}} + O \left( p^{- A} \right) ,
\end{equation} 
for any small $\epsilon > 0$ and any large $A > 0$, where
$$S(N)= \sum_{n \in \mathbb{Z}} \lambda_{f}(n)\lambda_{g}(n)V\left(\frac{n}{N}\right),$$
with $V$ a smooth bump function supported on $[\frac{1}{2},3]$ satisfying $V(x)=1$ for all $x\in [1,2]$.

To proceed further, in the delta method approach we separate oscillations of the Fourier coefficients $\lambda_{f}(n)$ and $\lambda_{g}(n)$, using the Jutila's circle method (see \cite{PA}, \cite{MZ}) similarly as done in the author's previous articles (see \cite{AG}, \cite{AG1}).

\subsection{Further set-up}
Let us apply the circle method directly to the smooth sum
$$S(N)= \mathop\sum_{n \in \mathbb{Z}}\lambda_{f}(n)\lambda_{g}(n)V\left(\frac{n}{N}\right),$$
 where the function $V$ is smooth, supported in $[1,2]$ with $V^{(j)}(x)\ll_{j} 1$. Now we will approximate the above sum $S(N)$ using Jutila's circle method (see \cite{PA}, \cite{MZ}) by the following sum :
\begin{equation}
\begin{split}
\tilde{S}(N)& =\frac{1}{L p}\mathop\sum_{q \in \Phi} \hspace{0.2cm}\sum_{b \bmod p} \, \, \sideset{}{^*}\sum_{a\bmod q}\mathop{\sum\sum}_{n,m \in \mathbb{Z}}\lambda_{f}(m)\lambda_{g}(n)e\left(\frac{a(n-m)}{p q}\right)\\
& \times e\left(\frac{b(n-m)}{p }\right)F(n,m) ,
\end{split}
\end{equation}
 where 
 $e_{q}(x)= e\left( \frac{x}{q}\right) = e^{2\pi i x /q}$, and 
$$F(x,y)= V\left(\frac{x}{N}\right)U\left(\frac{y}{N}\right)\frac{1}{2\delta}\int_{-\delta}^{\delta}e\left(\alpha \frac{(n-m)}{p}\right)d\alpha .$$
 Here $U$ is another smooth function having compact support in $(0,\infty)$, with $U(x)=1$ for $x$ in the support of $U$. Also we choose $\delta = \frac{p}{N}$ so that we have 
 $$\frac{\partial^{i+j}}{\partial^{i}x \partial^{j}y}F(x,y)\ll_{i,j}\frac{1}{N^{i+j}}.$$

$\hspace{-0.4 cm}$Then we have the following lemma :

\begin{lemma}\label{3.1} Let $\Phi \subset [Q,2Q]$, with 
$$L=\sum_{q\in \Phi}\phi(q) \gg Q^{2-\epsilon},$$
 and $\delta = \frac{p}{N} $. Then we must have 
\begin{equation}\label{m}
    S(N)= \tilde{S}(N) + O_{f,\epsilon}\left(N \sqrt{\frac{Q^2}{\delta L^2}}\right).
\end{equation}
\end{lemma}

\begin{proof}
Consider 
$$G(x)=\sum_{b \bmod p}\mathop{\sum\sum}_{n,m \in \mathbb{Z}}\lambda_{f}(m)\lambda_{g}(n) V\Big( \frac{n}{N}\Big) U\Big( \frac{m}{N}\Big)e\left(x\frac{(n-m)}{p}\right)e\left(b\frac{(n-m)}{p}\right) .$$

One can see that $S(N)= \int_{0}^{1}G(x) dx$ and $\Tilde{S}(N)=\int_{0}^{1}\tilde{I}_{\Phi , \delta}(x)G(x) dx$, where $\widetilde{I}(\alpha ): = \frac{1}{2\delta L} \mathop\sum_{q\in \Phi}\sum_{d \bmod q }^{\star }\mathcal{I}_{d/q}(\alpha )$ 
 and $\mathcal{I}_{d/q}$ is the indicator function of the interval $[\frac{d}{q}-\delta ,\frac{d}{q} + \delta ]$, $Q:= N^{1/2 + \epsilon}$ and $L \asymp Q^{2-\epsilon}$.

\noindent
Hence
\begin{equation}
\begin{split}
\Big| S(N)-\Tilde{S}(N)\Big| & \leq \frac{1}{p}\sum_{b \bmod p}\int_{0}^{1}\Big| 1- \tilde{I}_{\Phi , \delta}(x) \Big| \Big| \sum_{n \in \mathbf{Z}} \lambda_{g}(n)e\left(\frac{(x+b)n}{p}\right) V\Big( \frac{n}{N} \Big) \Big| \\
& \times \Big| \sum_{m \in \mathbb{Z}} \lambda_{f}(m)  e\left(\frac{(-x-b)m}{p}\right)U\Big( \frac{m}{N}\Big)\Big| dx .
\end{split}
\end{equation}

\noindent
 For the $n$-sum in the middle we have the Wilton-type bound
$$\sum_{n \in \mathbb{Z}} \lambda_{g}(n)e\left(\frac{(x+b)n}{p}\right) V\Big( \frac{n}{N} \Big) \ll_{g,\epsilon} N^{\frac{1}{2}+\epsilon}.$$
 Using the Cauchy-Schwarz's inequality to the remaining sum we get
 \begin{equation}
 \begin{split}
 \Big| S(N)-\Tilde{S}(N)\Big|& \ll \frac{N^{\frac{1}{2}+\epsilon}}{p} \Big( \sum_{b \bmod p} \Big[ \int_{0}^{1}\Big| 1- \tilde{I}_{\Phi , \delta}(x) \Big|^2 dx\Big]^{1/2} \\
 & \hspace{3cm}\times \Big[ \int_{0}^{1} \Big| \sum_{m \in \mathbb{Z}} \lambda_{f}(m)   e\left(-\frac{(x+b)m}{p}\right)U\Big( \frac{m}{N}\Big)\Big|^2 dx\Big]^{1/2} \big).
 \end{split}
 \end{equation}

For the last $m$-sum we open the absolute value square and execute the integral along with the additive character sum. So we are left with only the diagonal, which has size $N$. For the other sum we use \eqref{2.2}. It follows that
\begin{equation}
\Big| S(N)-\Tilde{S}(N)\Big|\ll N \sqrt{\frac{Q^2}{\delta L^2}}.
\end{equation}
\end{proof}

\subsection{Error term}\label{et} As for subconvexity, we need to save $N$ with something more, i.e., we must have $\delta L^2 > Q^2$, i.e., $\frac{1}{Q} \gg \delta > \frac{1}{Q^2}$, i.e., $Q > \sqrt{\frac{N}{p}}$. Here we choose $Q = \sqrt{\frac{N}{p}}p^{\eta + \epsilon}=p^{\frac{3}{2}+\eta+\epsilon}$. Hence the error term in Lemma \eqref{3.1} is bounded by $O\left(\frac{N}{p^{\eta}}\right)$.

\section{Sketch of the proof}
 Here in the sketch we will suppress the weight function for notational simplicity and we will closely follow the method done in \cite{AG}, \cite{AG1}, \cite{Mun}. We write the sum as 
 $$\mathbf{S}= \mathop{\sum\sum}_{n,m \sim p^{4}}\lambda_{f}(m)\lambda_{g}(n)\delta_{n,m},$$
 where $\delta_{n,m}$ is the Kronecker $\delta$-symbol. Here to get an inbuilt bilinear structure in the circle method itself, we need to use a more flexible version of the circle method - the one investigated by Jutila (see \cite{PA}, \cite{MZ}). This version comes with an error term which is satisfactory, as we shall find out, as long as we allow the moduli to be slightly larger than $\sqrt{N}$. Up to an admissible error we see that $\mathbf{S}$ is given by
$$\mathbf{S}=\mathop{\sum\sum}_{n,m \sim N}\lambda_{f}(m) \lambda_{g}(n)\int_{\mathbb{R}}\tilde{I}(\alpha )e((n-m)\alpha)d\alpha ,$$

\noindent
 where $\widetilde{I}(\alpha ): = \frac{1}{2\delta L} \mathop\sum_{q\in \Phi}\sum_{d \bmod q}^{\star }\mathcal{I}_{d/q}(\alpha )$ 
 and $\mathcal{I}_{d/q}$ is the indicator function of the interval $[\frac{d}{q}-\delta ,\frac{d}{q} + \delta ]$, $Q:= N^{1/2 + \epsilon}$ and $L \asymp Q^{2-\epsilon}$.

 Trivial bound at this stage yields $N^{2+\epsilon}$ and we need to establish the bound $N^{1-\theta}$ for some $\theta >0$, i.e., roughly speaking we need to save $N+$ something more. Observe that by our choice of $Q$, there is no analytic oscillation in the weight function $e((n-m)\alpha )$. Hence their weights can be dropped in our sketch. Here also we are assuming that $p \nmid q$. At first glue the congruence classes modulo $p$ and modulo $q$ to get $p \phi(q)$ many congruence classes modulo $p q$. In this sketch for simplicity we assume that all of the congruence classes are reduced modulo $p q$. So instead of $S(N)$ we consider the sum
 $$\sum_{q \sim Q} \, \, \sideset{}{^*}\sum_{a\bmod p q}\sum_{m \sim N}\lambda_{f}(m)e\left(\frac{-am}{p q}\right)\sum_{n \sim N}\lambda_{g}(n)e\left(\frac{an}{p q}\right) .$$
The $\mathrm{GL}(2)$ Voronoi on the $n$-sum yields
$$\sum_{n \sim N}\lambda_{g}(n)e\left(\frac{an}{p q}\right) \rightsquigarrow  \frac{N}{pq}\sum_{n \sim \frac{( p Q)^{2}}{N}}\lambda_{g}(n)e\left(\frac{-\overline{a}n}{p q}\right) ,$$
with a savings $\frac{N}{p Q}$.

Now we need to use the $\mathrm{GL}(2)$ Voronoi for the $m$-sum in a tricky way as $p^2$ appears in the modulus, but the modulus is not divisible by $p^4$, the level of the form. So the standard Voronoi summation formulas available in the literature (see \cite{KMV}) does not work directly here. This, however, can be overcome easily using the work of Atkin-Li (see \cite{AL}). Munshi (see \cite{Mun4}) used a trick for the Voronoi summation formula and we are going to use that trick here. The main point is that as $f$ is $p$-primitive, the twisted form $f \times \chi = f_{\chi}$ has level $[p^4 , p ^2 ] = p^4$, for any character $\chi$ modulo $p^2$. Hence the length of the dual sum is given by $m \rightsquigarrow \frac{p^4 Q^2}{N}$. So here Voronoi yields 
$$\sum_{m \sim N}\lambda_{f}(m)e\left(-\frac{am}{p q}\right) \rightsquigarrow  \frac{N}{p^2 q}\sum_{m \sim \frac{p^4 Q^2}{N}}\lambda_{f}(m)e\left(\frac{\overline{a}\overline{p^{3}}m}{q}\right) \mathbf{C}_{1}(a,m,q; p ) ,$$
where $\mathbf{C}_{1}(a,m,q; p )$ is an incomplete character sum roughly of the form
$$\sideset{}{^*}\sum_{\psi \bmod p}\psi(m\overline{a} \, \overline{q})\eta(f\times \overline{\psi})g(\psi) ,$$
with the savings $\frac{N}{pQ}$ provided that this incomplete character sum has a square root cancellation and it comes with a loss of $\sqrt{p}$.

\noindent
Now we consider the $a \bmod q$ sum which reduces to
$$\sideset{}{^*}\sum_{a \bmod q}e\left(\frac{-\overline{a}\overline{p}n}{ q}\right) e\left(\frac{\overline{a}\overline{p^{3}}m}{q}\right) \rightsquigarrow \mathcal{I}_{m \equiv p^2 n \bmod q},$$
where $\mathcal{I}_{x}$ is the characteristic function of $x$. So here we save $Q$. Also executing the sum over $a \bmod p$, it reduces to 
$$\sideset{}{^*}\sum_{a\bmod p}e\left(\frac{-\overline{a}n\overline{q}}{p }\right)\sideset{}{^*}\sum_{\psi \bmod p}\psi(m\overline{a} \, \overline{q})\eta(f\times \overline{\psi})g(\psi)$$
$$\rightsquigarrow \sideset{}{^*}\sum_{\psi \bmod p}\psi(m\overline{a} \, \overline{q})\eta(f\times \overline{\psi})g(\psi)^{2} \rightsquigarrow p^2 S(m,-n,;p) .$$
Now assuming the square root cancellation in the last sum we see that so far we have saved
$$\frac{1}{\sqrt{p}}\times\frac{N}{p Q}\times \frac{N}{p^{2}Q} \times Q\sqrt{p} \times \sqrt{p} = \frac{N}{p^{\eta}} .$$

\noindent
So we are on the boundary and need to save $p^{\eta}$ and something more. Now we are using the Cauchy-Schwarz's inequality, change of the variables and the Poisson summation formula to get our result.

\section{Voronoi summation formulae}
At first we consider the case $p \nmid b$ (other case will be dominated by this case as in that case we will have a saving of size $p$) and also we assume that each member of $\Phi$ is coprime to $p$, the modulus of the character $\chi$.

\noindent
Now let us consider
\begin{equation}\label{5.1}
     \Tilde{S}_{\neq 0}(N)= \frac{1}{L p}\sum_{q \in \Phi} \, \,  \sideset{}{^*}\sum_{a \bmod p q}S(a,q,x,f) T(a,q,x,g) ,
    \end{equation}
 where 
 $$S(a,q,x,g) := \sum_{n =1}^{\infty}\lambda_{g}(n)e_{p q}(an) e\left( \frac{xn}{p }\right)V\left(\frac{n}{N}\right) ,$$
 and
 $$T(a,q,x,f) := \sum_{m =1}^{\infty}\lambda_{f}(m)e_{p q}(-am) e\left( -\frac{xm}{p }\right)U\left(\frac{m}{N}\right) .$$

\noindent
so that we have 
\begin{equation}\label{s2}
    \Tilde{S}(N)= \frac{1}{2\delta}\int_{-\delta}^{\delta}\Tilde{S}_{x}(N) \, \mathrm{d}x .
\end{equation}

\subsection{Voronoi for the $n$-sum}\label{vn}
\begin{lemma}\label{gv} 
We have
\begin{equation}\label{eq12}
S\left(a,q,x,g\right)= \frac{N^{3/4}}{\sqrt{p}q^{1/2}}\sum_{1\leq n\ll N_{0}}\frac{\lambda_{g}(n)}{n^{1/4}}e\left(-\frac{\overline{a}n}{p q}\right)\mathcal{I}(n,x,q;p ) + O(N^{-2024}),
\end{equation}
where $N_{0} := p^{1+2\eta}$  and $\mathcal{I}(n,x,q;p )$ is given by
\begin{equation}\label{i}
    \mathcal{I}(n,x,q;p ):= \int_{\mathbb{R}}V(y)e\left( \frac{Nxy}{p } \pm \frac{4\pi}{p q}\sqrt{Nny}\right) W_g \left(\frac{4\pi \sqrt{Nny}}{p q}\right) \, \mathrm{d}y ,
\end{equation}
where $W_{g}$ is a smooth nice function.

\end{lemma}
\begin{proof}
Applying the Voronoi summation formula \eqref{V} to the $n$-sum of the equation \eqref{5.1}, then we have
$$\sum_{n \in \mathbb{Z}}\lambda_{g}(n)e\left(\frac{an}{p q}\right)e\left( \frac{xn}{p }\right)V\left(\frac{n}{N}\right) = \frac{1}{p q}\sum_{n \in \mathbb{Z}}\lambda_{g}(n)e\left(-\frac{\overline{a}n}{p q}\right)$$
 $$\hspace{4.5cm}\times\int_{\mathbb{R}}V\left(\frac{y}{N}\right)e\left( \frac{xy}{p }\right)J_{k_{g}-1}\left(\frac{4\pi \sqrt{ny}}{p q}\right) \, \mathrm{d}y, $$
 where $J_{k_{g}-1}$ is the Bessel function. By changing $y \mapsto Ny$ and using the decomposition \ref{bessel}, 
 $$J_{k_g -1}(x)= \frac{W_g (x)}{\sqrt{x}}e(x)+ \frac{\bar{W_g }(x)}{\sqrt{x}}e(-x) ,$$
  where $W_g (x)$ is a nice function, the right hand side integral becomes
 $$\frac{N^{3/4}}{q^{1/2}\sqrt{p}} \int_{\mathbb{R}}V(y)e\left( \frac{Nxy}{p } \pm \frac{4\pi}{p q}\sqrt{Nny}\right) W_g \left(\frac{4\pi \sqrt{Nny}}{p q}\right) \, \mathrm{d}y . $$
By repeated integration by parts we see that, this integral is negligibly small if $|n|\gg \frac{p^{{2}}q^{2} N^\epsilon}{N} =: N_{q}$ with $N_{q}\asymp \frac{p^{2}Q^{2} N^\epsilon}{N}=p^{1+2\eta}=: N_{0}$. Hence the lemma follows. 
 \end{proof}

\subsection{Voronoi for the $m$-sum}\label{vm}
Now let us consider the sum over $m$,
$$T(a,q,x,f) = \sum_{m =1}^{\infty}\lambda_{f}(m)e_{p q}(-am) e\left( \frac{-xm}{p }\right)U\left(\frac{m}{N}\right) .$$
Here we considered $p \nmid q$. To apply the Voronoi summation formula, we need to use results from \cite{AL}. Now for $p \nmid u$,
$$e\left( \frac{u}{p}\right) = \frac{1}{\phi(p)}\sum_{\chi \bmod p}\, g(\chi)\overline{\chi}(u)= -\frac{1}{\phi( p)}+ \frac{1}{\phi(p)}\sideset{}{^*}\sum_{\chi \bmod p}\, g(\chi)\overline{\chi}(u) ,$$
where $g(\chi)$ stands for the Gauss sum for $\chi$. Now as $p\nmid q$ so writing $a= a q \overline{q}+ap \overline{p}$ and also as $\lambda_{f}(p)=0$, so noting $p\nmid m$ in $T$, we have $T= T_{1} + T_{2}$, where
$$T_{1} := -\frac{1}{\phi(p)}\sum_{n =1}^{\infty}\lambda_{f}(m)e\left(-\frac{am\overline{p}}{q}\right) e\left( \frac{xm}{p }\right)U\left(\frac{m}{N}\right) ,$$
and
\begin{equation}\label{as2}
    T_{2} := \frac{1}{\phi(p)}\sideset{}{^*}\sum_{\chi \bmod p}\, g(\chi)\chi(q\overline{a})\sum_{m =1}^{\infty}\lambda_{f}(m)e\left(-\frac{am\overline{p}}{q}\right)e\left( \frac{xm}{p}\right)U\left(\frac{m}{N}\right) .
\end{equation}

\noindent
Here note that the contribution of the sum $T_{1}$ is dominated (as it comes with a prior saving of $p$) by the contribution of $T_{2}$. So we will examine the contribution of $T_{2}$.
\begin{lemma}\label{mf}
   We have
    \begin{equation}
    \begin{split}
        T_{2} &= \frac{2\pi i^{k_{f}}N}{q p^2 \phi(p)}\sideset{}{^*}\sum_{\chi \bmod p}\, g(\chi)\chi(\overline{qa})\eta(f\times \overline{\chi})\sum_{m=1}^{\infty}\lambda_{f}(m)\chi(m)e\left(\frac{\overline{a p^{3}} \, m}{q}\right)\\ 
        & \times\int_{0}^{\infty}U(y)e\left(\frac{xNy}{p}\right)J_{k_{f}}\left(\frac{4\pi\sqrt{mNy}}{qp^{2}}\right) + \, \mathrm{O}(N^{-2024}).
    \end{split}
    \end{equation} 
\end{lemma}
\begin{proof}
At first consider the sum
\begin{equation}\label{vf}
    \sum_{m=1}^{\infty}\lambda_{f}(m)\overline{\chi}(m)e\left(-\frac{a\overline{p}m}{q }\right)w(m)=\sum_{m=1}^{\infty}\lambda_{h}(m)e\left(-\frac{a\overline{p}m}{q }\right)w(m) \, \mathrm{d}y,
\end{equation}
where $h = f\times \overline{\chi}$. Now by proposition $3.1$ of \cite{AL}, we must have $h \in M_{k}(M, \overline{\chi}^{2})$ where $M= \lcm\left(\textit{level of }f , \textit{ conductor of } \chi \right)$ so for our case we have $M= \lcm( p^{4} , p^{2}) = p^{4}$. Also by the Theorem $3.1$ and Proposition $4.1$ of \cite{AL}, $h$ is a newform of level $p^{4}$ and nebentypus $\overline{\chi}^2$. Then by the Voronoi summation formula given in \cite{KMV}, the sum \eqref{vf} transforms into,
$$2\pi i^{k_{f}}\frac{N\overline{\chi}^{2}(q)\eta(h)}{qp^{2}}\sum_{m=1}^{\infty}\overline{\lambda_{h}(m)}e\left(\frac{\overline{a p^{3}}m}{q}\right) \int_{0}^{\infty}w(y)J_{k_{f}-1}\left(\frac{4\pi\sqrt{my}}{qp^{2}}\right) .$$
Substituting this in \eqref{vf} we get our lemma.
\end{proof}

Now extracting the oscillation of the Bessel function using \ref{bessel} as done for the $n$-sum and noting that $\phi(p )\asymp p$, we note that the $m$-sum $T_{2}$ in the Lemma \eqref{mf} is given by a sum of two similar terms of the form
\begin{equation}\label{fv}
\begin{split}
    \frac{N^{3/4}}{q^{1/2}p^{2}}\sideset{}{^*}\sum_{\chi \bmod p}\, g(\chi)\chi(\overline{qa})\eta(f\times \overline{\chi})\sum_{m=1}^{\infty}\frac{\lambda_{f}(m)\chi(m)}{m^{1/4}}e\left(\frac{\overline{a p^{3}} \, m}{q}\right) \mathcal{J}(m,x,q;p ),
 \end{split}
\end{equation}
where 
\begin{equation}\label{j}
    \mathcal{J}(m,x,q;p ) = \int_{0}^{\infty}U(y)e\left(-\frac{Nxy}{p^{2}}\pm \frac{2\sqrt{mNy}}{p^{2}q}\right)\, \mathrm{d}y .
\end{equation}
By repeated integration by parts we see that, this integral is negligibly small if $|m| \gg \frac{p^{4} q^2}{N}:= M_{q}\asymp \frac{p^{4} Q^2}{N}=p^{3+2\eta}=: M_{0}$ and we will consider the trivial bound for $\mathcal{J}(m,x,q;p )$ whenever needed.
\section{Further estimation}
Substituting \eqref{gv}, \eqref{fv} in \eqref{5.1} we arrive at
\begin{equation}\label{sn}
     \Tilde{S}_{\neq 0}(N)= \frac{N^{3/2}}{L p^\frac{7}{2}}\sum_{q \in \Phi} \frac{1}{q}  \sum_{1\leq n\ll N_{0}}\frac{\lambda_{g}(n)}{n^{1/4}}\sum_{m\ll M_{0}}\frac{\lambda_{f}(m)}{m^{1/4}}\mathcal{C}\mathcal{I}(n,x,q;p )\mathcal{J}(m,x,q;p ) ,
\end{equation}
where the pseudo-character sum $\mathcal{C}$ is given by
\begin{equation}\label{ch}
    \mathcal{C}(n,m,p,q):= \sideset{}{^*}\sum_{a \bmod p q}e\left(-\frac{\overline{a}n}{p q}\right)e\left(\frac{\overline{a p^{3}} \, m}{q}\right) \sideset{}{^*}\sum_{\chi \bmod p}\, g(\chi)\chi(\overline{qa}m)\eta(f\times \overline{\chi}) .
\end{equation}

\begin{lemma}\label{cl}
    We have
    \begin{equation}\label{chpq}
        \mathcal{C}(n,m,p,q)= \mathcal{C}_{p}(m,n;p )\, \, \mathcal{C}_{q} ,
    \end{equation}
    where 
    \begin{equation}\label{chp}
        \mathcal{C}_{p}(m,n;p ) = \phi(p)S(m,-n;p ) ,
    \end{equation}
    and 
    \begin{equation}\label{chq}
         \mathcal{C}_{q} =   \sum_{\substack{d \big| q \\ m \equiv p^{2}n \bmod d}}d \mu \left( \frac{q}{d} \right) ,
    \end{equation}

    \noindent
    and $\gcd (mn,p)=1$.
\end{lemma}
\begin{proof}
At first splitting the character sum in \eqref{ch} we get that
\begin{equation}\label{cp2}
    \mathcal{C}_{p}(m,n;p) = \sideset{}{^*}\sum_{a \bmod p }e\left(-\frac{\overline{a}\overline{q}n}{p }\right)\sideset{}{^*}\sum_{\chi \bmod p}\, g(\chi)\chi(\overline{qa}m)\eta(f\times \overline{\chi}) , 
\end{equation}
and 
\begin{equation}\label{cq2}
    \mathcal{C}_{q} = \sideset{}{^*}\sum_{a \bmod q }e\left(-\frac{\overline{ap}n}{ q}\right)e\left(\frac{\overline{a p^{{3}}} \, m}{q}\right) .
\end{equation}
Now the last character sum reduces to
\begin{equation}
\begin{split}
    \mathcal{C}_{q} & = \sideset{}{^*}\sum_{a \bmod q }e\left(-\frac{\overline{a}\overline{p}n}{ q}\right)e\left(\frac{\overline{a p^{3}} \, m}{q}\right) \\
    & = \sideset{}{^*}\sum_{a \bmod q }e\left(\overline{ap^{3}} \, \frac{(- n p^{2}+  \, m)}{q}\right) \\
    & = \sum_{\substack{d \big| q \\ m \equiv  p^{2} n \bmod d}}d \mu \left( \frac{q}{d} \right) .
\end{split}
\end{equation}
 Now consider the incomplete character sum modulo $p$. Changing variables and pushing the sum over $a$ inside and using the equation \eqref{ep} we have
\begin{equation*}
\begin{split}
     \mathcal{C}_{p}(m,n;p) & = \sideset{}{^*}\sum_{\chi \bmod p}\, g(\chi)\eta(f\times \overline{\chi})\sideset{}{^*}\sum_{a \bmod p }\chi(\overline{qa}m)e\left(-\frac{\overline{a}\overline{q}n}{p }\right) \\
     & = \sideset{}{^*}\sum_{\chi \bmod p}\, g(\chi)^{2}\eta(f\times \overline{\chi})\chi(m\overline{n}) \\
     & = \varepsilon(\pi_{f})\sideset{}{^*}\sum_{\chi \bmod p}g(\chi)^{2} \, \chi(m\overline{n}) \\
     & = \varepsilon(\pi_{f})\sideset{}{^*}\sum_{\chi \bmod p} \, \chi(m\overline{n}) \, \left( \sum_{a \bmod p}\chi(a)e\left(\frac{a}{p} \right)\right)^2 \\
     & = \varepsilon(\pi_{f})\mathop{\sum\sum}_{a_{1},a_{2}\bmod p} e\left(\frac{a_{1}}{p} \right)e\left(-\frac{a_{2}}{p} \right)\,\sideset{}{^*}\sum_{\chi \bmod p} \chi(m\overline{n}) \chi (a_{1})\chi (a_{2}) \\
     & = \varepsilon(\pi_{f})\phi(p)S(m,-n;p) - \varepsilon(\pi_{f}) \mathcal{I}_{m \equiv n \bmod p}.
\end{split}
\end{equation*}
\noindent
Here one can see that the sum \eqref{sn} with the second term from the previous line comes with an extra savings of $p$. So we will go with the first term only. Also here we note that $\gcd (mn,p)=1$.
\end{proof}

\section{Further computation}

 \subsection{Cauchy-Schwarz inequality} For the time being, we will focus on the case when $r=1$. At first we choose the set of moduli $\Phi$ to be the product set $\Phi_{1}\Phi_{2}$, where $\Phi_{i}$ consists of primes in the dyadic segment $[Q_{i}, 2Q_{i}]$ (and coprime to $p$) for $i = 1, 2$ and $Q_1 Q_2 = Q = \sqrt{\frac{N}{p}} p^{\eta }$. Also, we pick $Q_1$ and $Q_2$ (whose optimal sizes will be determined later) so that the collections $\Phi_1$ and $\Phi_{2}$ are disjoint. Now let us take $m= p^2 n + dr$ , so that $\gcd (r,p)=1$ as $\gcd (qm,p)=1$.

\noindent
We have

\begin{equation}\label{dq}
      \begin{split}
      \Tilde{S}_{\neq 0}(N):=\Sigma_{q}  & : = \frac{ N^{3/2}}{L p^\frac{5 }{2}}\sum_{q\in \Phi}\frac{1}{q} \sum_{d| q=q_{1}q_{2}}d\mu\left( \frac{q}{d}\right)\sum_{1\leq n\ll N_0}\left(\lambda_{g}(n)\frac{\mathcal{J}(n,x,q;p )}{n^{1/4}}\right)\\
         & \times \sum_{\substack{1\leq r \ll \frac{M_{0}}{d} \\ \gcd (r,p)=1}}\frac{\lambda_{f}( p^2 n +dr)\mathcal{I}(p^2 n+dr,x,q;p ) }{(p^2 n+dr)^{1/4}}S(dr,-n;p) 
     \end{split}
\end{equation}

\noindent
Several cases occur.

\subsection{Case $1$} Let $d= q $ (or $d=1$).

\begin{equation}
      \begin{split}
     \Sigma_{q} & =  \frac{ N^{3/2}}{L p^\frac{5 }{2}}\sum_{q\in \Phi} \sum_{1\leq n\ll N_0}\left(\lambda_{g}(n)\frac{\mathcal{J}(n,x,q;p )}{n^{1/4}}\right)\\
         & \times \sum_{\substack{1\leq r \ll \frac{M_{0}}{Q}\\ \gcd (rq,p)=1}}\frac{\lambda_{f}( p^2 n +qr)\mathcal{I}(p^2 n+qr,x,q;p ) }{(p^2 n+qr)^{1/4}}S(qr,-n;p) \\
         & \ll  \frac{ N^{3/2}}{L p^\frac{5 }{2}}\sum_{q\in \Phi}\sum_{\substack{1\leq r \ll \frac{M_{0}}{Q}\\ \gcd (rq,p)=1}}\Big| \\
         & \times \sum_{1\leq n\ll N_0}\left(\lambda_{g}(n)\frac{\mathcal{J}(n,x,q;p )}{n^{1/4}}\right)\frac{\lambda_{f}( p^2 n +qr)\mathcal{I}(p^2 n+qr,x,q;p ) }{(p^2 n+qr)^{1/4}}S(qr,-n;p) \Big|\\
         & \ll \frac{ N^{3/2}}{L p^\frac{5 }{2}}\sum_{\substack{1\leq r \ll M_{0}\\ \gcd(r,p)=1}}\Big| \\
         & \times \sum_{1\leq n\ll N_0}\lambda_{g}(n)\frac{\mathcal{J}(n,x,q;p )}{n^{1/4}}\frac{\lambda_{f}( p^2 n +r)\mathcal{I}(p^2 n+r,x,q;p ) }{(p^2 n+r)^{1/4}}S(r,-n;p) \Big|,
     \end{split}
\end{equation}

\noindent
where in the last line we have attached the variables $qr\mapsto r$ where we already had $\gcd(rq,p)=1$.

\noindent
After using the Cauchy-Schwarz's inequality we have

\begin{equation}\label{dq1}
      \begin{split}
      \Sigma_{q}   & \ll \frac{ N^{3/2}M_{0}^{1/2}}{L p^\frac{5 }{2}} \Big(\sum_{\substack{1\leq r \ll M_{0} \\ \gcd (r,p)=1}}\Big|\sum_{1\leq n \ll N_{0}}\lambda_{g}(n)\frac{\mathcal{J}(n,x,q;p )}{n^{1/4}}\\
      & \times \frac{\lambda_{f}( p^2 n +r) }{(p^2 n+r)^{1/4}} \mathcal{I}(p^2 n+r,x,q;p )S(r,-n;p)\Big|^{2}\Big)^{1/2}.
     \end{split}
\end{equation}

\noindent
Now consider the sum

\begin{equation}
\begin{split}
    \Sigma_{1}:& =\sum_{\substack{1\leq r \ll M_{0}\\ \gcd (r,p)=1}}\Big|\sum_{1\leq n \ll N_{0}}\lambda_{g}(n)\frac{\mathcal{J}(n,x,q;p )}{n^{1/4}}\frac{\lambda_{f}( p^2 n +r) }{(p^2 n+r)^{1/4}} \mathcal{I}(p^2 n+r,x,q;p )S(r,-n;p)\Big|^{2} \\
    & = \sum_{\substack{1\leq r \ll M_{0}\\ \gcd (r,p)=1}}\sum_{1\leq n \ll N_{0}}\lambda_{g}(n)\frac{\mathcal{J}(n,x,q;p )}{n^{1/4}}\frac{\lambda_{f}( p^2 n +r) }{(p^2 n+r)^{1/4}} \mathcal{I}(p^2 n+r,x,q;p )S(r,-n;p)\\
    & \times \sum_{1\leq n^\prime \ll N_{0}}\lambda_{g}(n^\prime )\frac{\mathcal{J}(n^\prime,x,q;p )}{n^{1/4}}\frac{\lambda_{f}( p^2 n^\prime +r) }{(p^2 n^\prime +r)^{1/4}} \mathcal{I}(p^2 n^\prime +r,x,q;p )S(r,-n^\prime ;p) .
\end{split}
\end{equation}

\noindent
We have two cases. 

\noindent
\emph{\textbf{Diagonal case.}} For the diagonal case $n= n^\prime$ we have

\begin{equation}\label{d1}
\begin{split}
    \Sigma_{1} =& \sum_{\substack{1\leq r \ll M_{0}\\ \gcd (r,p)=1}}\sum_{1\leq n \ll N_{0}}|\lambda_{g}(n)|^2 \frac{|\mathcal{J}(n,x,q;p )|^2}{n^{1/2}}\frac{|\lambda_{f}( p^2 n +r) |^2}{(p^2 n+r)^{1/2}} |\mathcal{I}(p^2 n+r,x,q;p )|^2 \, |S(r,-n;p)|^2 \\
    & \ll M_{0}^{1/2}N_{0}^{1/2}p = p^{3+2\eta},
\end{split}
\end{equation}

\noindent
where we have used the partial summation formula and the Ramanujan bound on average \eqref{RA}.

\noindent
\emph{\textbf{Off-diagonal case.}} For the off-diagonal case $n\neq n^\prime$ we have

\begin{equation}\label{sigma1}
\begin{split}
    \Sigma_{1} =& \sum_{\substack{1\leq r \ll M_{0}\\ \gcd (r,p)=1}}\sum_{1\leq n \ll N_{0}}\lambda_{g}(n)\frac{\mathcal{J}(n,x,q;p )}{n^{1/4}}\frac{\lambda_{f}( p^2 n +r) }{(p^2 n+r)^{1/4}} \mathcal{I}(p^2 n+r,x,q;p )S(r,-n;p)\\
    & \times \sum_{1\leq n^\prime \ll N_{0}}\lambda_{g}(n^\prime )\frac{\mathcal{J}(n^\prime ,x,q;p )}{{n^\prime}^{1/4}}\frac{\lambda_{f}( p^2 n^\prime +r) }{(p^2 n^\prime +r)^{1/4}} \mathcal{I}(p^2 n^\prime +r,x,q;p )S(r,-n^\prime ;p) .
\end{split}
\end{equation}

Using the Cauchy-Schwarz's inequality the above inequality becomes dominated by

\begin{equation}
\begin{split}
    \Sigma_{1} \ll & \sqrt{N_{0}}\Big(\sum_{1\leq n \ll N_{0}}\sum_{1\leq n^\prime \ll N_{0}}\Big|\sum_{\substack{1\leq r \ll M_{0}\\ \gcd (r,p)=1}}\frac{\lambda_{f}( p^2 n +r) }{(p^2 n+r)^{1/4}} \mathcal{I}(p^2 n+r,x,q;p )S(r,-n;p)\\
    & \times \frac{\lambda_{f}( p^2 n^\prime +r) }{(p^2 n^\prime +r)^{1/4}} \mathcal{I}(p^2 n^\prime +r,x,q;p )S(r,-n^\prime ;p) \Big|^2 \Big)^{1/2}.
\end{split}
\end{equation}

\noindent
Changing the variables $p^2 n+r \mapsto u$ and $n^\prime = n+\ell$ so that $\gcd (r,p)=1 \implies \gcd (u,p)=1$ and hence we have

\begin{equation}
\begin{split}
    \Sigma_{1} \ll & \sqrt{N_{0}}\Big(\sum_{1\leq n \ll N_{0}}\sum_{1\leq \ell \ll N_{0}}\Big|\sum_{\substack{1\leq u \ll M_{0}\\ \gcd (u,p)=1}}\frac{\lambda_{f}( u) }{u^{1/4}} \mathcal{I}(u,x,q;p )S(u,-n;p)\\
    & \times \frac{\lambda_{f}( u+\ell p^2 ) }{(u+ \ell p^2)^{1/4}} \mathcal{I}(u+\ell p^2,x,q;p )S(u,-n-\ell ;p) \Big|^2 \Big)^{1/2}.
\end{split}
\end{equation}

\noindent
Now opening the absolute value square we have

\begin{equation}
\begin{split}
    \Sigma_{2}  & =  \sum_{1\leq n \ll N_{0}}\sum_{1\leq \ell \ll N_{0}}\sum_{\substack{1\leq u \ll M_{0}\\ \gcd (u,p)=1}}\frac{\lambda_{f}( u) }{u^{1/4}} \mathcal{I}(u,x,q;p )S(u,-n;p)\\
    & \times \frac{\lambda_{f}( u+\ell p^2 ) }{(u+ \ell p^2)^{1/4}} \mathcal{I}(u+\ell p^2,x,q;p )S(u,-n-\ell ;p)  \\
    & \times \sum_{\substack{1\leq u^\prime \ll M_{0}\\ \gcd (u^\prime ,p)=1}}\frac{\lambda_{f}( u^\prime ) }{{u^\prime}^{1/4}} \mathcal{I}(u^\prime ,x,q;p )S(u^\prime ,-n;p)\\
    & \times \frac{\lambda_{f}( u^\prime +\ell p^2 ) }{(u^\prime + \ell p^2)^{1/4}} \mathcal{I}(u^\prime +\ell p^2,x,q;p )S(u^\prime ,-n-\ell ;p) .
\end{split}
\end{equation}

\noindent
We have several cases:

\noindent
\emph{\textbf{Diagonal case.}} For $u \equiv u^\prime \bmod p$, 

\begin{equation}
\begin{split}
    \Sigma_{2}  & \ll  \sum_{1\leq n \ll N_{0}}\sum_{1\leq \ell \ll N_{0}}\sum_{\substack{1\leq u \ll M_{0}\\ \gcd (u,p)=1}}\frac{|\lambda_{f}( u)|^2 }{u^{1/2}} |\mathcal{I}(u,x,q;p )|^2 |S(u,-n;p)|^2\\
    & \times \frac{|\lambda_{f}( u+\ell p^2 )|^2 }{(u+ \ell p^2)^{1/2}} |\mathcal{I}(u+\ell p^2,x,q;p )|^2 |S(u,-n-\ell ;p)|^2 \\
    & \ll N_{0}^2 p^2 ,
\end{split}
\end{equation}

\noindent
where we have used the partial summation formula and the Ramanujan bound on average \eqref{RA}.

\noindent
So for this case, we have

\begin{equation}\label{d2}
\begin{split}
    \Sigma_{1} \ll & \sqrt{N_{0}}\Big( N_{0}^2 p^2\Big)^{1/2} \ll N_{0}^{3/2} p = p^{\frac{5}{2}+3\eta},
\end{split}
\end{equation}

\noindent
where we have used the partial summation formula and the Ramanujan bound on average \eqref{RA}.

\noindent
\emph{\textbf{Off-diagonal case.}} For $u \not\equiv u^\prime \bmod p$, 

\begin{equation}
\begin{split}
    \Sigma_{2}  & =  \sum_{1\leq \ell \ll N_{0}}\sum_{\substack{1\leq u \ll M_{0}\\ \gcd (u,p)=1}}\frac{\lambda_{f}( u) }{u^{1/4}} \mathcal{I}(u,x,q;p )\\
    & \times \frac{\lambda_{f}( u+\ell p^2 ) }{(u+ \ell p^2)^{1/4}} \mathcal{I}(u+\ell p^2,x,q;p )  \\
    & \times \sum_{\substack{1\leq u^\prime \ll M_{0}\\ \gcd (u^\prime ,p)=1}}\frac{\lambda_{f}( u^\prime ) }{{u^\prime}^{1/4}} \mathcal{I}(u^\prime ,x,q;p )\\
    & \times \frac{\lambda_{f}( u^\prime +\ell p^2 ) }{(u^\prime + \ell p^2)^{1/4}} \mathcal{I}(u^\prime +\ell p^2,x,q;p ) \\
    & \times \left( \sum_{1\leq n\ll N_{0}}S(u,-n;p) S(u,-n-\ell ;p)S(u^\prime ,-n;p)S(u^\prime ,-n-\ell ;p)\right) .
\end{split}
\end{equation}

Now taking $n = pm+\beta$ where $\beta \bmod p$ and then we do the Poisson summation formula on the $n$-sum and note that as the size of the conductor is $p \ll N_{0}= p^{1+2\eta}$ so the zero frequency will remain only, i.e., the case for $m=0$ and so we have

\begin{equation}\label{sigma21}
\begin{split}
    \Sigma_{2}  & =\sum_{1\leq \ell \ll N_{0}}\sum_{\substack{1\leq u \ll M_{0}\\ \gcd (u,p)=1}}\frac{\lambda_{f}( u) }{u^{1/4}} \mathcal{I}(u,x,q;p )\\
    & \times \frac{\lambda_{f}( u+\ell p^2 ) }{(u+ \ell p^2)^{1/4}} \mathcal{I}(u+\ell p^2,x,q;p )  \\
    & \times \sum_{\substack{1\leq u^\prime \ll M_{0}\\ \gcd (u^\prime ,p)=1}}\frac{\lambda_{f}( u^\prime ) }{{u^\prime}^{1/4}} \mathcal{I}(u^\prime ,x,q;p )\\
    & \times \frac{\lambda_{f}( u^\prime +\ell p^2 ) }{(u^\prime + \ell p^2)^{1/4}} \mathcal{I}(u^\prime +\ell p^2,x,q;p ) \\
    & \times \frac{N_{0}}{p}\left( \sum_{\beta \bmod p}S(u,-\beta;p) S(u,-\beta-\ell ;p)S(u^\prime ,-\beta;p)S(u^\prime ,-\beta-\ell ;p)\right) .
\end{split}
\end{equation}

\noindent
Now consider the character sum

\begin{equation}\label{a}
    \begin{split}
       \mathcal{A}:& = \sum_{\beta \bmod p}S(u,-\beta;p) S(u,-\beta-\ell ;p)S(u^\prime ,-\beta;p)S(u^\prime ,-\beta-\ell ;p)\\
      & = \sideset{}{^*}\sum_{a_{1}\bmod p}e\left( \frac{\overline{a_{1}}u}{p}\right)\sideset{}{^*}\sum_{a_{2}\bmod p}e\left( \frac{\overline{a_{2}}u-a_{2}\ell }{p}\right)\sideset{}{^*}\sum_{a_{3}\bmod p}e\left( \frac{\overline{a_{3}}u^\prime }{p}\right)\sideset{}{^*}\sum_{a_{4}\bmod p}e\left( \frac{\overline{a_{4}}u^\prime -a_{4}\ell }{p}\right) \\
      & \times \sum_{\beta \bmod p}e\left( \frac{-\beta a_{1}}{p}\right)e\left( \frac{-\beta a_{2}}{p}\right) e\left( \frac{-\beta a_{3}}{p}\right)e\left( \frac{-\beta a_{4}}{p}\right).
    \end{split}
\end{equation}

\noindent
Now from the $\beta$ sum we have

\begin{equation*}
    \begin{split}
      & a_{1}+a_{2}+a_{3}+a_{4}\equiv 0 \bmod p \\
      & a_{1}+a_{2}+a_{3}\equiv -a_{4} \bmod p .
    \end{split}
\end{equation*}

\noindent
Hence the character sum becomes

\begin{equation*}
    \begin{split}
       \mathcal{A} & = p\sideset{}{^*}\sum_{a_{1}\bmod p}\, \, \sideset{}{^*}\sum_{a_{2}\bmod p}\, \, \sideset{}{^*}\sum_{a_{3}\bmod p}\, \, e\left( \frac{\overline{a_{1}}u+\overline{a_{2}}u+\overline{a_{3}}u^\prime  -(\overline{a_{1}+a_{2}+a_{3}})u^\prime +a_{1}\ell+a_{3}\ell}{p}\right). 
    \end{split}
\end{equation*}

\noindent
Now several cases occur.

\noindent
\emph{\textbf{ Case $1$.}} For the case $\ell = 0$, the above sum becomes  

\begin{equation*}
    \begin{split}
      \mathcal{A}=& p \sideset{}{^*}\sum_{a_{2}\bmod p}\, \, \sideset{}{^*}\sum_{a_{3}\bmod p}\, \, \sideset{}{^*}\sum_{a_{1}\bmod p}\, \, e\left( \frac{\overline{a_{1}}(u+\overline{a_{2}}u+\overline{a_{3}}u)  -\overline{a_{1}}(\overline{1+a_{2}+a_{3}})u}{p}\right) \\
      & = p \sideset{}{^*}\sum_{a_{2}\bmod p}\, \, \sideset{}{^*}\sum_{a_{3}\bmod p} \left(1-p\mathbb{I}_{u+\overline{a_{2}}u+\overline{a_{3}}u-(\overline{1+a_{2}+a_{3}})u\equiv 0 \bmod p} \right)\\
      & = p\left(p^2 - p\mathbb{I}_{u \equiv 0 \bmod p}-p\mathbb{I}_{1+\overline{a_{2}}+\overline{a_{3}}-(\overline{1+a_{2}+a_{3}})\equiv 0 \bmod p} \right).
    \end{split}
\end{equation*}

\noindent
Now as the case $u \equiv 0 \bmod p$ is not possible so we will focus mainly on the last equation in which $a_{3}$ can be determined by $a_{2}$ so we have 

\begin{equation*}
    \begin{split}
      & p \sideset{}{^*}\sum_{a_{2}\bmod p}\, \, \sideset{}{^*}\sum_{a_{3}\bmod p}\, \, \sideset{}{^*}\sum_{a_{1}\bmod p}\, \, e\left( \frac{\overline{a_{1}}(u+\overline{a_{2}}u+\overline{a_{3}}u)  -\overline{a_{1}}(\overline{1+a_{2}+a_{3}})u}{p}\right) \\
      & = p\left(p^2 - p\mathbb{I}_{u \equiv 0 \bmod p}-p(p-1) \right)\\
      & = p\left( - p\mathbb{I}_{u \equiv 0 \bmod p}+p \right) \ll p^2 .
    \end{split}
\end{equation*}

\noindent
Hence from \eqref{sigma21} and using the partial summation formula along with the Ramanujan bound on average \eqref{RA} we have

    \begin{equation*}
\begin{split}
    \Sigma_{2}  & \ll M_{0} \times \frac{N_{0}}{p}\times p^{2}  \ll p^{5+4\eta}
\end{split}
\end{equation*}

\noindent
and so from \eqref{sigma1} we have 

\begin{equation}\label{o2}
    \Sigma_{1}\ll \sqrt{N_{0}}\times p^{\frac{5}{2}+2\eta} = p^{3+3\eta}
\end{equation}

\noindent
\emph{\textbf{ Case $2$.}} For the case $\ell \equiv  0 \bmod p$ and $\ell \neq 0$, changing the variables $\ell \mapsto p\ell$, the character sum becomes  

\begin{equation*}
    \begin{split}
     \mathcal{A}= & p \sideset{}{^*}\sum_{a_{2}\bmod p}\, \, \sideset{}{^*}\sum_{a_{3}\bmod p}\, \, \sideset{}{^*}\sum_{a_{1}\bmod p}\, \, e\left( \frac{\overline{a_{1}}(u+\overline{a_{2}}u+\overline{a_{3}}u^\prime)  -\overline{a_{1}}(\overline{1+a_{2}+a_{3}})u^\prime}{p}\right) \\
      & = p \sideset{}{^*}\sum_{a_{2}\bmod p}\, \, \sideset{}{^*}\sum_{a_{3}\bmod p} \left(1-p\mathbb{I}_{u+\overline{a_{2}}u+\overline{a_{3}}u-(\overline{1+a_{2}+a_{3}})u\equiv 0 \bmod p} \right)\\
      & = p\left(p^2 - p\mathbb{I}_{u \equiv 0 \bmod p}-p\mathbb{I}_{1+\overline{a_{2}}+\overline{a_{3}}-(\overline{1+a_{2}+a_{3}})\equiv 0 \bmod p} \right).
    \end{split}
\end{equation*}

\noindent
Now for the case $u \equiv 0 \bmod p$ we will save $p$ so we will focus mainly on the last equation in which $a_{3}$ can be determined by $a_{2}$ so we have 

\begin{equation*}
    \begin{split}
      & p \sideset{}{^*}\sum_{a_{2}\bmod p}\, \, \sideset{}{^*}\sum_{a_{3}\bmod p}\, \, \sideset{}{^*}\sum_{a_{1}\bmod p}\, \, e\left( \frac{\overline{a_{1}}(u+\overline{a_{2}}u+\overline{a_{3}}u)  -\overline{a_{1}}(\overline{1+a_{2}+a_{3}})u}{p}\right) \\
      & = p\left(p^2 - p\mathbb{I}_{u \equiv 0 \bmod p}-p(p-1) \right)\\
      & = p\left( - p\mathbb{I}_{u \equiv 0 \bmod p}+p \right) \ll p^2 .
    \end{split}
\end{equation*}

\noindent
Hence from \eqref{sigma21} and using the partial summation formula along with the Ramanujan bound on average \eqref{RA} we have

    \begin{equation*}
\begin{split}
    \Sigma_{2}  & \ll \frac{N_{0}}{p}\times M_{0} \times \frac{N_{0}}{p}\times p^{2}  \ll p^{5+6\eta}
\end{split}
\end{equation*}

\noindent
and so from \eqref{sigma1} we have 

\begin{equation}\label{o3}
    \Sigma_{1}\ll \sqrt{N_{0}}\times p^{\frac{5}{2}+3\eta} = p^{3+4\eta}
\end{equation}

\noindent
\emph{\textbf{ Case $3$.}} Now consider the case when $\ell \not\equiv 0 \bmod p$. For this case, the character sum can be written as follows.

\begin{equation*}
    \begin{split}
       \mathcal{A}&= \sum_{\beta \bmod p}S(u,-\beta;p) S(u,-\beta-\ell ;p)S(u^\prime ,-\beta;p)S(u^\prime ,-\beta-\ell ;p)\\
      & = p^2 \sum_{\beta \bmod p}\frac{S(1,-u\beta;p) }{\sqrt{p}}\frac{S(1,-u\beta-u\ell ;p) }{\sqrt{p}}\frac{S(1 ,-u^\prime\beta;p) }{\sqrt{p}}\frac{S(1 ,-u^\prime\beta-u^\prime\ell ;p)}{\sqrt{p}}\\
      &= p^2 \sum_{0 \leq n < p}\cos \left(\theta_{p,f_{1}(n)}\right) \, \cos \left(\theta_{p,f_{2}(n)}\right) \, \cos \left(\theta_{p,f_{3}(n)}\right) \, \cos \left(\theta_{p,f_{4}(n)}\right) ,
    \end{split}
\end{equation*}

\noindent
where 

$$f_{1}(n)=-un, \, f_{2}(n)=-un- u\ell , \, f_{3}(n)=-u^\prime n, \, f_{4}(n)=-u^\prime n -u^\prime\ell .$$

\noindent
Here we note that $u \not\equiv u^\prime \bmod p, \, u\not\equiv 0 \bmod p, \, u^\prime\not\equiv 0 \bmod p, \, \ell \not\equiv 0 \bmod p$. Then for $r=4, \, k_{i}=1$ we have from \cite[Lemma $2.1$]{FMRS}

\begin{equation}\label{c4}
    \mathcal{A} \ll C(4) p^{5/2},
\end{equation}

\noindent
for some non-zero constant $C(4)$. So we need to take $p>C(4)$.

\noindent
Hence, using the partial summation formula along with the Ramanujan bound on average \eqref{RA} we have 

\begin{equation}\label{o1}
    \Sigma_{1} \ll p^{3+4\eta} + p^{\frac{15}{4}+3\eta}.
\end{equation}

\noindent
So, from \eqref{d1} and \eqref{o1} we have 

\begin{equation}\label{m2}
    \Sigma_{1} \ll p^{3+4\eta} + p^{\frac{15}{4}+3\eta} +p^{\frac{7}{2}+3\eta}.
\end{equation}

\noindent
Hence, from \eqref{dq1} we have

\begin{equation}\label{dq2}
      \begin{split}
      \Sigma_{q}   & \ll \frac{ N^{3/2}M_{0}^{1/2}}{L p^\frac{5 }{2}} \Big(p^{3+4\eta} + p^{\frac{15}{4}+3\eta} +p^{\frac{5}{2}+3\eta}\Big)^{1/2}\\
      & \ll \sqrt{N }\times \left(p^{\frac{3}{2}+\eta}+p^{\frac{15}{8}+\frac{\eta}{2}}+ p^{\frac{5}{4}+\frac{\eta}{2}}\right).
     \end{split}
\end{equation}

\subsection{Case $2$} Let $d= q_{1} $ (or $d= q_{2} $).

\begin{equation}
      \begin{split}
     \Sigma_{q_{1}} & = = \frac{ N^{3/2}}{L p^\frac{5 }{2}}\sum_{q\in \Phi} \frac{1}{q_{2}}\sum_{1\leq n\ll N_0}\left(\lambda_{g}(n)\frac{\mathcal{J}(n,x,q;p )}{n^{1/4}}\right)\\
         & \times \sum_{1\leq r \ll \frac{M_{0}}{Q_{1}}}\frac{\lambda_{f}( p^2 n +q_{1}r)\mathcal{I}(p^2 n+q_{1}r,x,q;p ) }{(p^2 n+q_{1}r)^{1/4}}S(q_{1}r,-n;p) \\
         & \ll  \frac{ N^{3/2}}{L p^\frac{5 }{2}}\sum_{q\in \Phi}\frac{1}{q_{2}}\sum_{1\leq r \ll \frac{M_{0}}{Q_{1}}}\Big| \\
         & \times \sum_{1\leq n\ll N_0}\left(\lambda_{g}(n)\frac{\mathcal{J}(n,x,q;p )}{n^{1/4}}\right)\frac{\lambda_{f}( p^2 n +q_{1}r)\mathcal{I}(p^2 n+q_{1}r,x,q;p ) }{(p^2 n+q_{1}r)^{1/4}}S(q_{1}r,-n;p) \Big|\\
         & \ll \frac{ N^{3/2}}{L p^\frac{5 }{2}}\sum_{q_{2}\in \Phi_{2}}\frac{1}{q_{2}}\sum_{Q_{1}\leq r \ll M_{0}}\Big| \\
         & \times \sum_{1\leq n\ll N_0}\lambda_{g}(n)\frac{\mathcal{J}(n,x,q;p )}{n^{1/4}}\frac{\lambda_{f}( p^2 n +r)\mathcal{I}(p^2 n+r,x,q;p ) }{(p^2 n+r)^{1/4}}S(r,-n;p) \Big|,
     \end{split}
\end{equation}

\noindent
where in the last line we have attached the variables $q_{1}r\mapsto r$. 

\noindent
Now processing similarly as done for the above case we have that

\begin{equation}\label{dq3}
      \begin{split}
      \Sigma_{q_{1}}   & \ll \frac{ N^{3/2}M_{0}^{1/2}}{L p^\frac{5 }{2}} \left(p^{\frac{3}{2}+\eta}+p^{\frac{15}{8}+\frac{\eta}{2}}+ p^{\frac{5}{4}+\frac{\eta}{2}}\right).
     \end{split}
\end{equation}

\section{Further computation}

From \eqref{dq}, \eqref{dq2}, \eqref{dq3} and the error term we have that

\begin{equation*}
\begin{split}
     \Tilde{S}_{\neq 0}(N) &\ll \sqrt{N }\times \left( p^{\frac{3}{2}+\eta}+p^{\frac{15}{8}+\frac{\eta}{2}}+ p^{\frac{5}{4}+\frac{\eta}{2}} + \frac{p^2}{p^\eta}\right)\\
     & \ll \sqrt{N }\times \left( p^{\frac{3}{2}+\eta}+p^{\frac{15}{8}+\frac{\eta}{2}}+ \frac{p^2}{p^\eta}\right).
    \end{split}
   \end{equation*}

\noindent
Now compairing the last two terms we have $\eta = \frac{1}{12}$ and so 

\begin{equation*}
\begin{split}
     \Tilde{S}_{\neq 0}(N) &\ll \sqrt{N }\times p^{\frac{23}{12}}.
    \end{split}
   \end{equation*}

\noindent
This gives the Theorem \ref{mt}.

\section{Details for Maass forms}

Above we have considered that $f, g$ are holomorphic cusp forms. Now instead of considering them as holomorphic cusp forms, if we take $f$ (and /or $g$) as Maass cusp form, then we need to go through some minor adjustments (starting from preliminaries in section \ref{sec2}) in the following sense:

\begin{itemize}
    \item We need to adjust the gamma factors in the functional equation \ref{V} accordingly (see \cite[Chap. $5.11$]{IK}) and also we need to adjust the weight functions with appropriate weight functions accordingly.

    \item We need to adjust the Voronoi summation formula stated in the Lemma \ref{V} accordingly (see \cite[Theorem $A.4$]{KMV}). Also for the application of the Voronoi summation formula (see subsections \ref{vn} and \ref{vm}), we need to adjust the Bessel functions accordingly with $\pm$ sign introduced on each of the $n$ and $m$-sums. Then we will continue with the rest of the arguments as \ref{bessel} also holds for the other Bessel functions with appropriate weight functions $W_{k}$. 
.
    \item The Ramanujan bound on average \eqref{RA} that we have used for the Fourier coefficients $\lambda_{f}(n)$ and $\lambda_{g}(m)$ will be unchanged.

\end{itemize}

\
{}

\end{document}